\documentclass[12pt]{amsart}
\usepackage{fullpage,url,amssymb}
\usepackage[all]{xy} 

\DeclareFontEncoding{OT2}{}{} 
\newcommand{\textcyr}[1]{%
 {\fontencoding{OT2}\fontfamily{wncyr}\fontseries{m}\fontshape{n}\selectfont #1}}
\newcommand{\Sha}{{\mbox{\textcyr{Sh}}}}

\usepackage{color}

\newcommand{\defi}[1]{\textsf{#1}} 

\newcommand{\Aff}{{\mathbb A}}

\newcommand{\G}{{\mathbb G}}

\newcommand{\PP}{{\mathbb P}}
\newcommand{\Q}{{\mathbb Q}}

\newcommand{\Z}{{\mathbb Z}}

\newcommand{\kbar}{{\overline{k}}}

\newcommand{\Bbar}{{\overline{B}}}
\newcommand{\Cbar}{{\overline{C}}}
\newcommand{\Vbar}{{\overline{V}}}
\newcommand{\Xbar}{{\overline{X}}}

\newcommand{\Zbar}{{\overline{Z}}}

\newcommand{\Adeles}{{\mathbb A}}

\newcommand{\calC}{{\mathcal C}}

\newcommand{\calE}{{\mathcal E}}

\newcommand{\calL}{{\mathcal L}}
\newcommand{\calM}{{\mathcal M}}

\newcommand{\calV}{{\mathcal V}}

\newcommand{\LL}{{\mathcal L}}
\newcommand{\OO}{{\mathcal O}}

\DeclareMathOperator{\supp}{supp}

\DeclareMathOperator{\Hom}{Hom}

\DeclareMathOperator{\Gal}{Gal}
\DeclareMathOperator{\Ind}{Ind}

\DeclareMathOperator{\Br}{Br}

\DeclareMathOperator{\Sym}{Sym}
\DeclareMathOperator{\Div}{Div}
\DeclareMathOperator{\Pic}{Pic}

\DeclareMathOperator{\Spec}{Spec}

\newcommand{\Et}{\operatorname{\bf Et}}

\newcommand{\descent}{{\operatorname{descent}}}
\newcommand{\connected}{{\operatorname{connected}}}
\newcommand{\conts}{{\operatorname{conts}}}

\newcommand{\et}{{\operatorname{et}}}

\newcommand{\injects}{\hookrightarrow}
\newcommand{\isom}{\simeq}

\newcommand{\Intersection}{\bigcap} 
\newcommand{\Union}{\bigcup} 
\newcommand{\tensor}{\otimes}
\newcommand{\directsum}{\oplus} 
\newcommand{\Directsum}{\bigoplus} 

\newcommand{\isomto}{\overset{\sim}{\rightarrow}}

\newtheorem{theorem}{Theorem}[section]
\newtheorem{lemma}[theorem]{Lemma}

\newtheorem{proposition}[theorem]{Proposition}

\theoremstyle{definition}

\newtheorem{question}[theorem]{Question}

\theoremstyle{remark}
\newtheorem{remark}[theorem]{Remark}

\usepackage[
        backref,
        pdfauthor={Bjorn Poonen}, 
]{hyperref}
\usepackage[alphabetic,backrefs,lite]{amsrefs} 

\begin{document}

\title[Brauer-Manin obstruction applied to \'etale covers]{Insufficiency of the Brauer-Manin obstruction applied to \'etale covers}
\author{Bjorn Poonen}
\thanks{This research was supported by NSF grant DMS-0301280.}
\address{Department of Mathematics, University of California, 
        Berkeley, CA 94720-3840, USA}
\email{poonen@math.berkeley.edu}
\urladdr{http://math.berkeley.edu/\~{}poonen/}
\date{June 8, 2008}

\begin{abstract}
Let $k$ be any global field of characteristic not~$2$.
We construct a $k$-variety $X$ such that $X(k)$ is empty, 
but for which the emptiness cannot be explained 
by the Brauer-Manin obstruction
or even by the Brauer-Manin obstruction applied to finite \'etale covers.
\end{abstract}

\maketitle

\section{Introduction}\label{S:introduction}

\subsection{Background}

Call a variety \defi{nice} if it is smooth, projective, 
and geometrically integral.
(See Sections \ref{S:notation} and~\ref{S:obstructions} for 
further terminology used here.)
Let $X$ be a nice variety over a global field $k$.
If $X$ has a $k$-point, then $X$ has a $k_v$-point for every place
$v$ of $k$; i.e., the set $X(\Adeles)$ of adelic points is nonempty.
The converse, known as the \defi{Hasse principle}, does not always hold,
as has been known at least since the 1940s: it can fail for genus-$1$ 
curves, for instance \cites{Lind1940,Reichardt1942}.
Manin \cite{Manin1971} showed that the Brauer group of $X$
can often explain failures of the Hasse principle:
one can define a subset
$X(\Adeles)^{\Br}$ of $X(\Adeles)$ that contains $X(k)$,
and $X(\Adeles)^{\Br}$ can be empty even when $X(\Adeles)$ is nonempty.

Conditional results \cites{Sarnak-Wang1995,Poonen-complete-intersection2001}
predicted that this Brauer-Manin obstruction
was insufficient to explain all failures of 
the Hasse principle.
But the insufficiency was proved only in 1999, 
when a groundbreaking paper of Skorobogatov~\cite{Skorobogatov1999}
constructed a variety for which one could prove 
$X(\Adeles)^{\Br} \ne \emptyset$ and $X(k)=\emptyset$.
He showed that for a bielliptic surface $X$,
the set $X(\Adeles)^{\et,\Br}$
obtained by applying the Brauer-Manin obstruction to finite \'etale covers
of $X$ could be empty even when $X(\Adeles)^{\Br}$ was not.

\subsection{Our result}

We give a construction to show that even this combination of
finite \'etale descent and the Brauer-Manin obstruction 
is insufficient to explain all failures of the Hasse principle.
Our argument does not use \cite{Skorobogatov1999},
so it also gives a new approach to constructing
varieties for which the Brauer-Manin obstruction is insufficient
to explain the failure of the Hasse principle.

The idea behind our construction can be described in a few lines,
though the details will occupy the rest of the paper.
Start with a nice curve $C$ such that $C(k)$ is finite and nonempty.
Construct a nice $k$-variety $X$ 
with a morphism $\beta \colon X \to C$ such that
\begin{enumerate}
\item[(i)] For each $c \in C(k)$, 
the fiber $X_c:=\beta^{-1}(c)$ violates the Hasse principle.
\item[(ii)] Every finite \'etale cover of $X$ 
arises from an \'etale cover of $C$.
\item[(iii)] The map $\beta$ induces an isomorphism $\Br C \isomto \Br X$,
and this remains true after base extension by any finite \'etale 
morphism $C' \to C$.
\end{enumerate}
Properties (ii) and~(iii) imply that $X(\Adeles)^{\et,\Br}$
is the inverse image under $\beta$ of $C(\Adeles)^{\et,\Br}$,
which contains the nonempty set $C(k)$.
Then by~(i), we have $X(\Adeles)^{\et,\Br} \ne \emptyset$
but $X(k)=\emptyset$.

Our $X$ will be a $3$-fold, and the general fiber of $\beta$ 
will be a Ch\^atelet surface (a kind of conic bundle over $\PP^1$).

\subsection{Commentary}

Suppose in addition that the Jacobian $J$ of $C$ is such that
the Mordell-Weil group $J(k)$ and the Shafarevich-Tate group $\Sha(J)$
are both finite.
(For instance, these hypotheses are known to hold if $C$ is any
elliptic curve over $\Q$ of analytic rank $0$.)
Then $C(\Adeles)^{\Br}$ is essentially (ignoring some technicalities
regarding the connected components at archimedean places)
equal to $C(k)$.
(Scharaschkin and Skorobogatov independently observed that this
follows from the comparison of the Cassels-Tate pairing with the
Brauer evaluation pairing in \cite{Manin1971}: 
see \cite{Skorobogatov2001}*{\S6.2} for related results,
and \cite{Stoll2007}*{Theorem~8.6} for a significant generalization.)
Thus $X(\Adeles)^{\Br}$ is essentially a subset of
$\Union_{c \in C(k)} X_c(\Adeles)$.

Also, $X_c(\Adeles)^{\Br}=\emptyset$ for each $c \in C(k)$
(all failures of the Hasse principle for Ch\^atelet surfaces are explained by 
the Brauer-Manin obstruction \cites{CT-Sansuc-SD1987I,CT-Sansuc-SD1987II}).
But the elements of $\Br X_c$ used to obstruct $k$-points on the fiber $X_c$
do not extend to elements of $\Br X$,
so it does not follow that $X(\Adeles)^{\Br}$ is empty.

\subsection{Outline of the paper}

Section~\ref{S:notation} introduces some basic notation.
Section~\ref{S:obstructions} recalls some cohomological obstructions
to rational points, and discusses how they relate to one another.
Our $X$, a Ch\^atelet surface bundle over $C$, 
will be constructed as a conic bundle over $C \times \PP^1$;
Section~\ref{S:conic bundles} describes the type of conic bundle we need,
and Sections \ref{S:etale covers} and \ref{S:Brauer group}
compute the \'etale covers and Brauer group of this conic bundle.
The Brauer group calculations involve some group cohomology lemmas,
which have been relegated to an appendix.
Section~\ref{S:construction} constructs the particular $X$,
and Sections \ref{S:no BM} and~\ref{S:no descent}
compute $X(\Adeles)^{\Br}$ and $X(\Adeles)^{\et,\Br}$,
respectively.

\section{Notation}\label{S:notation}

Given a field $k$, 
we fix a separable closure $\kbar$ of $k$ and define $G_k:=\Gal(\kbar/k)$.
For any $k$-variety $V$, define $\Vbar:=V \times_k \kbar$.
For any nice variety $V$, let $\kappa(V)$ be the function field.
If $D$ is a divisor on a nice variety $V$, 
let $[D]$ be its class in $\Pic V$.

An \defi{algebraic group} over $k$ 
is a smooth group scheme of finite type over $k$.
If $G$ is an algebraic group over $k$ and $X$ is a $k$-variety,
let $H^1(X,G)$ be the cohomology set 
defined using \v{C}ech $1$-cocycles for the \'etale topology;
it parameterizes isomorphism classes of torsors over $X$ under $G$.
If moreover $G$ is commutative,
then for any $i \in \Z_{\ge 0}$ 
define $H^i(X,G)$ as the usual \'etale cohomology group;
this is compatible with the previous definition when $i=1$.
Let $\Br X$ be the cohomological Brauer group $H^2(X,\G_m)$.

By a \defi{global field} we mean either a finite extension of $\Q$
or the function field of a nice curve over a finite field.
If $k$ is a global field, let $\Adeles$ be its ad\`ele ring.

\section{Cohomological obstructions to rational points}\label{S:obstructions}

Let $k$ be a global field.
Let $X$ be a nice $k$-variety.

\subsection{Brauer-Manin obstruction}
(See \cite{Skorobogatov2001}*{\S5.2}.) 
There is an evaluation pairing
\[
	\Br X \times X(\Adeles) \to \Q/\Z,
\]
and $X(\Adeles)^{\Br}$ is defined as the set of elements of $X(\Adeles)$
that pair with every element of $\Br X$ to give $0$.
The reciprocity law for $\Br k$ 
implies $X(k) \subseteq X(\Adeles)^{\Br}$.
In particular, if $X(\Adeles) \ne \emptyset$
but $X(\Adeles)^{\Br} = \emptyset$,
then $X$ violates the Hasse principle.

\subsection{Descent obstruction}

(See \cite{Skorobogatov2001}*{\S5.3}.)
If $G$ is a (not necessarily connected) linear algebraic group over $k$,
and $f\colon Y \to X$ is a right torsor under $G$,
then any $1$-cocycle $\sigma \in Z^1(k,G)$
gives rise to a ``twisted'' right torsor $f^\sigma \colon Y^\sigma \to X$
under a twisted form $G^\sigma$ of $G$.
Moreover, the isomorphism type of the torsor depends only on
the cohomology class $[\sigma] \in H^1(k,G)$.
It is not hard to show that
\begin{equation}
\label{E:coverings}
        X(k) = \Union_{[\sigma] \in H^1(k,G)} f^\sigma(Y^\sigma(k)).
\end{equation}
Therefore $X(k)$ is contained in the set
\[
        X(\Adeles)^f := 
        \Union_{[\sigma] \in H^1(k,G)} f^\sigma(Y^\sigma(\Adeles)).
\]
Define
\[
        X(\Adeles)^\descent := \Intersection X(\Adeles)^f
\]
where the intersection is taken over all linear algebraic groups $G$
and all right torsors $Y \to X$ under $G$.
If $X(\Adeles)^\descent =\emptyset$, then we say that 
there is a \defi{descent obstruction} to the existence of a rational point.

\subsection{Brauer-Manin obstruction applied to \'etale covers}

For reasons that will be clearer in Section~\ref{S:comparisons},
it is interesting to combine 
descent for torsors under finite \'etale group schemes
with the Brauer-Manin obstruction:
define
\[
	X(\Adeles)^{\et,\Br}:= \Intersection_{\substack{\text{$G$ finite}\\ f \colon Y \to X}} \quad \Union_{[\sigma] \in H^1(k,G)} f^\sigma(Y^\sigma(\Adeles)^{\Br}),
\]
where the intersection is taken over all finite \'etale group schemes $G$
over $k$ and all right torsors $f \colon Y \to X$ under $G$.
We have $X(k) \subseteq X(\Adeles)^{\et,\Br} \subseteq X(\Adeles)^{\Br}$,
where the first inclusion follows from~\eqref{E:coverings},
and the second follows from taking $G=\{1\}$ and $Y=X$ in the definition
of $X(\Adeles)^{\et,\Br}$.

\subsection{Comparisons}
\label{S:comparisons}

Let $X(\Adeles)^\connected$ 
be defined in the same way as $X(\Adeles)^\descent$,
but using only {\em connected} linear algebraic groups
instead of all linear algebraic groups.
Harari \cite{Harari2002}*{Th\'eor\`eme 2(2)}
showed that $X(\Adeles)^{\Br} \subseteq X(\Adeles)^\connected$,
at least if $k$ is a number field.
In other words, the Brauer-Manin obstruction is strong enough to subsume
all descent obstructions from connected linear algebraic groups.
Also, an arbitrary linear algebraic group 
is an extension of a finite \'etale group scheme
by a connected linear algebraic group, so one might ask:
\begin{question}
\label{Q:descent vs et-Br}
Does $X(\Adeles)^{\et,\Br} \subseteq X(\Adeles)^{\descent}$ hold
for every nice variety $X$ over a number field?
\end{question}

This does not seem to follow formally from Harari's result.
But Demarche, in response to an early draft of this paper,
has announced a positive answer~\cite{Demarche2008-preprint}.

\section{Conic bundles}\label{S:conic bundles}

In this section, $k$ is any field of characteristic not $2$.
Let $B$ be a nice $k$-variety.
Let $\calL$ be a line sheaf on $B$.
Let $\calE$ be the rank-3 vector sheaf $\OO \directsum \OO \directsum \LL$
on $B$.
Let $a \in k^\times$
and let $s \in \Gamma(B,\calL^{\tensor 2})$ be a nonzero global section.
The zero locus of 
\[
        1 \directsum (-a) \directsum (-s) 
        \in \Gamma(B,\OO \directsum \OO \directsum \LL^{\tensor 2}) 
        \subseteq \Gamma(B,\Sym^2 \calE)
\]
in $\PP\calE$ is a projective geometrically integral scheme $X$
with a morphism $\alpha\colon X \to B$.
If $U$ is a dense open subscheme of $B$ with a trivialization
$\LL|_U \isom \OO_U$
and we identify $s|_U$ with an element of $\Gamma(U,\OO_U)$,
then the affine scheme defined by $y^2-az^2=s|_U$ in $\Aff^2_U$
is a dense open subscheme of $X$.
Therefore we call $X$ the \defi{conic bundle given by $y^2-az^2=s$}.
In the special case where $B=\PP^1$, $\LL = \OO(2)$,
and the homogeneous form $s \in \Gamma(\PP^1,\OO(4))$ is separable,
$X$ is called the \defi{Ch\^atelet surface given by $y^2-az^2 =s(x)$},
where $s(x) \in k[x]$ denotes a dehomogenization of $s$.

Returning to the general case, we let $Z$ be the subscheme $s=0$ of $B$.
Call $Z$ the \defi{degeneracy locus} of the conic bundle.
Each fiber of $\alpha$ above a point of $B-Z$ is a smooth plane conic,
and each fiber above a geometric point of $Z$ 
is a union of two projective lines crossing transversely at a point.
A local calculation shows that if $Z$ is smooth over $k$, 
then $X$ is smooth over $k$.

\begin{lemma}
\label{L:generic fiber}
The generic fiber $\Xbar_\eta$ of $\Xbar \to \Bbar$
is isomorphic to $\PP^1_{\kappa(\Bbar)}$.
\end{lemma}

\begin{proof}
It is a smooth plane conic,
and it has a rational point since $a$ is 
a square in $\kbar \subseteq \kappa(\Bbar)$.
\end{proof}

\section{\'Etale covers of conic bundles}\label{S:etale covers}

Given a variety $X$, let $\Et(X)$ be the category of finite \'etale 
covers of $X$.
We recall some well-known properties of $\Et(X)$:

\begin{lemma}
\label{L:Et}
Let $X$ be a regular integral variety over a field $k$.
Let $Z \subseteq X$ be a closed subvariety not equal to $X$ itself.
Let $U:=X-Z$.
Then
\begin{enumerate}
\item[(i)] 
The inclusion $U \injects X$ induces 
a fully faithful functor $\Et(X) \to \Et(U)$.
\item[(ii)] 
If the codimension of $Z$ in $X$ is at least $2$,
then $\Et(X) \to \Et(U)$ is an equivalence of categories.
\end{enumerate}
\end{lemma}

\begin{proof}\hfill
\begin{enumerate}
\item[(i)] Use the description of \'etale covers of a normal
integral scheme given in \cite{SGA1}*{I.10.3}.
\item[(ii)] This is ``purity of the branch locus'' 
\cite{SGA2}*{X.3.3 and X.3.4}. 
\qedhere
\end{enumerate}
\end{proof}

\begin{lemma}
\label{L:Et and P^1}
For any $k$-variety $X$, the projection $X \times \PP^1 \to X$
induces an equivalence of categories $\Et(X) \to \Et(X \times \PP^1)$.
\end{lemma}

\begin{proof}
This follows from the fact that $\PP^1$ has no finite \'etale covers
except for those induced by finite \'etale covers of $\Spec k$.
\end{proof}

\begin{lemma}
\label{L:etale covers of birational varieties}
A birational map $\phi\colon X \dashrightarrow Y$ between nice varieties
induces an equivalence of categories $\Et(Y) \to \Et(X)$.
\end{lemma}

\begin{proof}
Let $U = X-Z$ be the dense open subscheme of $X$ on which $\phi$ is defined.
Since $X$ is regular and $Y$ is proper, $Z$ has codimension at least $2$
in $X$.
We have functors $\Et(Y) \to \Et(U) \leftarrow \Et(X)$,
and the second is an equivalence of categories by Lemma~\ref{L:Et}(ii).
In the same way, $\phi^{-1}$ induces a functor $\Et(X) \to \Et(Y)$.
To see that these are inverse equivalences of categories,
use Lemma~\ref{L:Et}(i) to
compare with the equivalence $\Et(Y_1) \to \Et(X_1)$
where $X_1 \subseteq X$ and $Y_1 \subseteq Y$ are nonempty open subschemes
such that $\phi$ restricts to an isomorphism $X_1 \to Y_1$.
\end{proof}

\begin{proposition}
\label{P:etale covers of conic bundle}
Let $\alpha\colon X \to B$ be a conic bundle given by $y^2-az^2=s$
as in Section~\ref{S:conic bundles}, with smooth degeneracy locus.
Then $\Et(B) \to \Et(X)$ is an equivalence of categories.
\end{proposition}

\begin{proof}
By Lemma~\ref{L:generic fiber}, the nice $\kbar$-varieties
$\Xbar$ and $\Bbar \times \overline{\PP^1}$ are birational.
Lemmas \ref{L:Et and P^1} and~\ref{L:etale covers of birational varieties}
gives equivalences
$\Et(\Bbar) \to \Et(\Bbar \times \overline{\PP^1}) \to \Et(\Xbar)$.
The composite equivalence 
$\Et(\Bbar) \to \Et(\Xbar)$
is induced by the base extension of a $k$-variety morphism,
so descent data for descending covers over $\kbar$ to covers over $k$
correspond under $\Et(\Bbar) \to \Et(\Xbar)$.
Thus $\Et(B) \to \Et(X)$ is an equivalence too.
\end{proof}

\section{Brauer group of conic bundles}\label{S:Brauer group}

The calculations of this section are similar to 
well-known calculations that have been done for conic bundles over $\PP^1$: 
see \cite{Skorobogatov2001}*{\S7.1}, for instance.

\begin{lemma}
\label{L:H^1 of Pic}
Let $X \to B$ be as in Section~\ref{S:conic bundles}.
If the degeneracy locus $Z$ is nice, then the homomorphism 
$H^1(k,\Pic \Bbar) \stackrel{\alpha^*}\to H^1(k,\Pic \Xbar)$ 
is an isomorphism.
\end{lemma}

\begin{proof}
We compute $\Pic \Xbar$ in the following paragraphs 
by constructing a commutative diagram of $G_k$-modules
\begin{equation}
\label{E:key}
\xymatrix{
0 \ar[r] & \Z \ar[r]^-{\lambda_1} \ar@{=}[d] & \tilde{\Z}^2 \ar[r]^-{\lambda_2} \ar[d]_{\tau_1} & \tilde{\Z}^2 \ar[r]^-{\lambda_3} \ar[d]_{\tau_2} & \Z \ar[r] \ar@{=}[d]_{\deg} & 0 \\
0 \ar[r] & \Z \ar[r]^-{\rho_1} & \Pic \Bbar \directsum \tilde{\Z}^2 \ar[r]^-{\rho_2} & \Pic \Xbar \ar[r]^-{\rho_3} & \Pic \Xbar_\eta \ar[r] & 0 \\
}
\end{equation}
with exact rows.

Let $\tilde{\Z}^2$ be the induced module $\Ind_{G_{k(\sqrt{a})}}^{G_k} \Z$:
as a group it is $\Z^2$, and an element $\sigma \in G_k$ acts on an element
of it either trivially or by interchanging the coordinates,
according to whether $\sigma$ fixes $\sqrt{a}$ or not.

Call a divisor of $\Xbar$ \defi{vertical} if it is supported on
prime divisors lying above prime divisors of $\Bbar$,
and \defi{horizontal} otherwise.
The fiber of $\alpha$ above the generic point of $\Zbar$
consists of two intersecting copies of $\PP^1_{\kappa(\Zbar)}$,
so $\alpha^{-1}(\Zbar)$ is a union of two prime divisors $F_1$ and $F_2$
of $\Xbar$.

Choose $L \in \Div B$ with $[L]=\LL$.
Since $Z$ is the zero locus of $s \in \Gamma(B,\LL^{\tensor 2})$,
the divisor $Z-2L$ is the divisor of some function $g \in \kappa(B)^\times$.
Let $U := B - \supp(L)$.
Then $X$ has an open subscheme $X'$ given by $y^2 - a = g t^2$
in the affine space $\Aff^2_U$ with coordinates $t$ and $y$.
The restrictions of $F_1$ and $F_2$ in $\Div X'$ 
are given by $y-\sqrt{a}=g=0$ and $y+\sqrt{a}=g=0$;
we may assume that the former is $F_1$.
The Zariski closures in $\Xbar$ of the divisors 
given by $y-\sqrt{a}=t=0$ and $y+\sqrt{a}=t=0$
are horizontal; call them $H_1$ and $H_2$.
We choose a function $f \in \kappa(\Xbar)^\times$ 
that on the generic fiber 
induces an isomorphism $\Xbar_\eta \to \PP^1_{\kappa(\Bbar)}$
(the usual parameterization of a conic);
explicitly, we take
\[
        f:= \frac{y-\sqrt{a}}{t} = \frac{g t}{y+\sqrt{a}}.
\]
A straightforward calculation shows that the divisor of $f$ on $\Xbar$ is
\begin{equation}
\label{E:divisor}
        (f) = H_1 - H_2 + F_1 - \alpha^* L.
\end{equation}

{\em Bottom row:}
Define $\rho_1$ by $\rho_1(1)=(-2\calL,(1,1))$.
Define $\rho_2(\calM,(m,n))=\alpha^*\calM + m[F_1] + n[F_2]$.
Let $\rho_3$ be restriction.
Each $\rho_i$ is $G_k$-equivariant.
Given a prime divisor $D$ on $\Xbar_\eta$, 
its Zariski closure in $\Xbar$ restricts to give $D$ on $\Xbar_\eta$,
so $\rho_3$ is surjective.
The kernel of $\rho_3$ 
is generated by the classes of vertical prime divisors of $\Xbar$;
in fact, there is exactly one above each prime divisor of $\Bbar$ 
except that above $\Zbar \in \Div \Bbar$ we have $F_1, F_2 \in \Div \Xbar$.
This proves exactness at $\Pic \Xbar$ of the bottom row.
Since $s \in \Gamma(B,\LL^{\tensor 2})$,
we have $[Z]=2\LL$, and $[F_1]+[F_2] = \alpha^* [Z] = 2 \alpha^* \LL$.
Also, a rational function on $\Xbar$ with vertical divisor
must be the pullback of a rational function on $\Bbar$;
The previous two sentences prove exactness 
at $\Pic \Bbar \directsum \tilde{\Z}^2$.
Injectivity of $\rho_1$ is trivial,
so this completes the proof that the bottom row of~\eqref{E:key}
is exact.

{\em Top row:}
Define
\begin{align*}
  \lambda_1(m) &= (m,m) \\
  \lambda_2(m,n) &= (n-m,m-n) \\
  \lambda_3(m,n) &= m+n.
\end{align*}
These maps are $G_k$-equivariant and 
they make the top row of~\eqref{E:key} exact.

{\em Vertical maps:}
By Lemma~\ref{L:generic fiber}, 
we have an isomorphism $\deg \colon \Pic \Xbar_\eta \isom \Z$
of $G_k$-modules; this defines the rightmost vertical map in~\eqref{E:key}.
Define
\begin{align*}
   \tau_1(m,n) &= (-(m+n)\LL,(m,n)) \\
   \tau_2(m,n) &= m [H_1] + n [H_2].
\end{align*}
These too are $G_k$-equivariant.

Commutativity of the first square is immediate from the definitions.
Commutativity of the second square follows from~\eqref{E:divisor}.
Commutativity of the third square follows since $H_1$ and $H_2$
each meet the generic fiber $\Xbar_\eta$ in a single $\kappa(\Bbar)$-rational
point.
This completes the construction of~\eqref{E:key}.

We now take cohomology by applying results of Appendix~\ref{A:cohomology}.
Because of the vertical isomorphisms at the left and right ends 
of~\eqref{E:key}, the two rows define the same class $\xi \in H^2(G_k,\Z)$.
We have $H^0(G_k,\tilde{\Z}^2) = \Z \cdot (1,1)$,
and Shapiro's lemma yields 
\begin{equation}
\label{E:Shapiro}
        H^1(G_k,\tilde{\Z}^2)=H^1(G_{k(\sqrt{a})},\Z)=0,
\end{equation}
so Lemma~\ref{L:nonzero H^2} implies $\xi \ne 0$.
We are almost ready to apply Lemma~\ref{L:2-extension}
to the bottom row of~\eqref{E:key},
but first we must check the splitting hypotheses.
After restricting from $G_k$ to $G_{k(\sqrt{a})}$,
the injection $\rho_1$ is split by the projection
$\Pic \Bbar \directsum \Z^2 \to \Z$ onto the last factor,
and the surjection $\rho_3$ is split by the map sending a positive generator
of $\Pic \Xbar_\eta$ to $[H_1] \in \Pic \Xbar$.
Now Lemma~\ref{L:2-extension} yields an isomorphism
\[
        H^1(G_k,\Pic \Bbar \directsum \tilde{\Z}^2) \to H^1(G_k,\Pic \Xbar)
\]
and the first group equals $H^1(G_k,\Pic \Bbar)$ by~\eqref{E:Shapiro}.
\end{proof}

\begin{lemma}
\label{L:ruled}
If $W$ and $Y$ are nice $k$-varieties,
and $W$ is birational to $Y \times \PP^1$,
then the homomorphism $\Br Y \to \Br W$ 
induced by the composition $W \dashrightarrow Y \times \PP^1 \to Y$
is an isomorphism.
\end{lemma}

\begin{proof}
Use the birational invariance of the Brauer group
and the isomorphism $\Br(Y \times \PP^1) \isom \Br Y$.
\end{proof}

\begin{lemma}
\label{L:Br Xbar}
Let $X \to B$ be as in Section~\ref{S:conic bundles}.
If $\Br \Bbar=0$, then $\Br \Xbar=0$.
\end{lemma}

\begin{proof}
Apply Lemmas \ref{L:generic fiber} and~\ref{L:ruled}.
\end{proof}

\begin{proposition}
\label{P:Br X}
Let $X \to B$ be as in Section~\ref{S:conic bundles}.
Suppose in addition that 
\begin{itemize}
\item $k$ is a global field (still of characteristic not $2$),
\item the degeneracy locus $Z$ is nice,
\item $\Br \Bbar=0$, and
\item $X(\Adeles) \ne \emptyset$.
\end{itemize}
Then $\alpha^*\colon \Br B \to \Br X$ is an isomorphism.
\end{proposition}

\begin{proof}
The Hochschild-Serre spectral sequence yields an exact sequence
\[
        \Br k 
        \to \ker\left(\Br X \to \Br \Xbar \right) 
        \to H^1(k,\Pic \Xbar) 
        \to H^3(k,\G_m).
\]
Since $\Br k \to \Directsum_v \Br k_v$ is injective 
and $X(\Adeles)\ne \emptyset$,
the homomorphism $\Br k \to \Br X$ is injective.
By Lemma~\ref{L:Br Xbar}, we have $\Br \Xbar=0$.
Finally, $H^3(k,\G_m)=0$.
Thus we obtain a short exact sequence, the second row of 
\[
\xymatrix{
0 \ar[r] & \Br k \ar[r] \ar@{=}[d] & \Br B \ar[r] \ar[d] & H^1(k,\Pic \Bbar) \ar[d] \ar[r] & 0 \\
0 \ar[r] & \Br k \ar[r] & \Br X \ar[r] & H^1(k,\Pic \Xbar) \ar[r] & 0. \\
}
\]
The first row is obtained in the same way, and the vertical maps
are induced by $\alpha$.
The result now follows from Lemma~\ref{L:H^1 of Pic}.
\end{proof}

\begin{remark}
\label{R:alternative proof}
In response to an earlier draft of this paper,
Colliot-Th\'el\`ene has found an alternative proof
of Proposition~\ref{P:Br X}:
see~\cite{Colliot-Thelene2008-preprint}*{Proposition~2.1}.
This proof, which is a little shorter and works in slightly greater
generality, compares $\Br X$ and $\Br B$ using residue maps
instead of going through $H^1(k,\Pic \Xbar)$
and $H^1(k,\Pic \Bbar)$.
\end{remark}

\section{Construction}\label{S:construction}

{}From now on, $k$ is a global field of characteristic not $2$.
Fix $a \in k^\times$,
and fix relatively prime separable 
degree-$4$ polynomials $P_\infty(x), P_0(x) \in k[x]$
such that the (nice) Ch\^atelet surface $\calV_\infty$ given by
\[
        y^2 - a z^2 = P_\infty(x)
\]
over $k$ satisfies $\calV_\infty(\Adeles)\ne \emptyset$ 
but $\calV_\infty(k)=\emptyset$.
(Such Ch\^atelet surfaces exist over any global field $k$ 
of characteristic not~$2$: 
see \cite{Poonen-chatelet-preprint}*{Proposition~5.1 and \S11}.
If $k=\Q$, one may use the original example from~\cite{Iskovskikh1971},
with $a:=-1$ and $P_\infty(x):=(x^2-2)(3-x^2)$.)
Let $\tilde{P}_\infty(w,x)$ and $\tilde{P}_0(w,x)$ be the homogenizations
of $P_\infty$ and $P_0$.
Define $\LL:=\OO(1,2)$ on $\PP^1 \times \PP^1$
and define
\[
        s_1:=u^2 \tilde{P}_\infty(w,x) + v^2 \tilde{P}_0(w,x) 
                \in \Gamma(\PP^1 \times \PP^1,\LL^{\tensor 2}),
\]
where the two copies of $\PP^1$ have homogeneous coordinates $(u,v)$
and $(w,x)$, respectively.
Let $Z_1 \subset \PP^1 \times \PP^1$ be the zero locus of $s_1$.
Let $F \subset \PP^1$ be the (finite) branch locus of the 
first projection $Z_1 \to \PP^1$.

Let $\alpha_1 \colon \calV \to \PP^1 \times \PP^1$ 
be the conic bundle given by $y^2-az^2=s_1$,
in the terminology of Section~4.
Composing $\alpha_1$ with the first projection $\PP^1 \times \PP^1 \to \PP^1$
yields a morphism $\beta_1 \colon \calV \to \PP^1$
whose fiber above $\infty:=(1:0)$ is the Ch\^atelet surface $\calV_\infty$
defined earlier.

Let $C$ be a nice curve over $k$ 
such that $C(k)$ is finite and nonempty.
Choose a dominant morphism $\gamma \colon C \to \PP^1$,
\'etale above $F$, such that $\gamma(C(k))=\{\infty\}$.

Define the fiber product $X:= \calV \times_{\PP^1} C$ 
and morphisms $\alpha$ and $\beta$ as in the diagram
\[
\xymatrix{
X \ar[d]_{\alpha} \ar@/^2.5pc/[dd]^>>>>>>{\beta} \ar[rr] && \calV \ar[d]_{\alpha_1} \ar@/^2.5pc/[dd]^>>>>>>{\beta_1} \\
C \times \PP^1 \ar[rr]^{(\gamma,1)} \ar[d]_{1^{\operatorname{st}}} && \PP^1 \times \PP^1 \ar[d]_{1^{\operatorname{st}}} \\
C \ar[rr]^\gamma && \PP^1. \\
}
\]
Each map labeled $1^{\operatorname{st}}$ is the first projection.
Define $B:=C \times \PP^1$ 
and $s :=(\gamma,1)^* s_1 \in \Gamma(B,(\gamma,1)^* \OO(2,4))$.
Thus $X \stackrel{\alpha}\to B$ can alternatively be described as
the conic bundle given by $y^2-az^2=s$.
Its degeneracy locus $Z$ is $(\gamma,1)^* Z_1 \subset B$.

\section{No Brauer-Manin obstruction}\label{S:no BM}

We continue with the notation of Section~\ref{S:construction}.

\begin{lemma}
\label{L:curve Z}
The curve $Z$ is nice.
\end{lemma}

\begin{proof}
Since $P_0(x)$ and $P_\infty(x)$ are separable and have no common factor,
a short calculation shows that $Z_1$ is smooth over $k$.
We have $Z = Z_1 \times_{\PP^1} C$.
Since $Z_1$ and $C$ are smooth over $k$ 
and the branch loci of $Z_1 \to \PP^1$ and $C \to \PP^1$ do not intersect,
$Z$ is smooth too.
Since $Z_1$ is ample on $\PP^1 \times \PP^1$ and $\gamma$ is finite,
$Z$ is ample on $C \times \PP^1$.
Therefore, $Z$ is geometrically connected 
by \cite{Hartshorne1977}*{Corollary~III.7.9}.
Since $Z$ is also smooth, it is geometrically integral.
\end{proof}

Lemma~\ref{L:curve Z} and the sentence before Lemma~\ref{L:generic fiber}
imply that the $3$-fold $X$ is nice.

\begin{theorem}
\label{T:no BM obstruction}
We have $X(k)=\emptyset$, 
but $X(\Adeles)^{\Br}$ contains $\calV_\infty(\Adeles) \times C(k)$
and hence is nonempty.
\end{theorem}

\begin{proof}
Since $\gamma(C(k))=\{\infty\}$ and $\calV_\infty(k)=\emptyset$,
we have $X(k)=\emptyset$.

We have $\Br \Bbar = \Br(\Cbar \times \overline{\PP^1}) = \Br \Cbar = 0$,
by Lemma~\ref{L:ruled} and \cite{Grothendieck-Brauer3}*{Corollaire~5.8}.
Also, $X(\Adeles)$ contains $\calV_\infty(\Adeles) \times C(k)$,
so $X(\Adeles) \ne \emptyset$.
Thus Proposition~\ref{P:Br X} implies that $\Br B \to \Br X$ is an isomorphism.
Composing with the isomorphism $\Br C \to \Br B$ of Lemma~\ref{L:ruled}
shows that $\beta^* \colon \Br C \to \Br X$ is an isomorphism.
Hence, if $\beta_\Adeles \colon X(\Adeles) \to C(\Adeles)$
is the map induced by $\beta$,
then 
\[
        X(\Adeles)^{\Br} 
        = \beta_\Adeles^{-1}(C(\Adeles)^{\Br}) 
        \supseteq \beta_{\Adeles}^{-1}(C(k)) 
        = \calV_\infty(\Adeles) \times C(k).
\]
\end{proof}

\section{No Brauer-Manin obstruction applied to \'etale covers}\label{S:no descent}

We continue with the notation of Section~\ref{S:construction};
in particular, $X$ is the nice $3$-fold defined there.

\begin{theorem}
\label{T:no descent obstruction}
The set $X(\Adeles)^{\et,\Br}$ contains 
$\calV_\infty(\Adeles) \times C(k)$
and hence is nonempty.
\end{theorem}

\begin{proof}
Lemma~\ref{L:Et and P^1} applied to $C$
and Proposition~\ref{P:etale covers of conic bundle} 
applied to $X \to B = C \times \PP^1$
yield equivalences of categories
$\Et(C) \to \Et(C \times \PP^1) \to \Et(X)$.

Let $G$ be a finite \'etale group scheme over $k$,
and let $f\colon Y \to X$ be a right torsor under $G$.
The equivalence $\Et(C) \to \Et(X)$ implies that $f$ arises from
a right torsor $h \colon \calC \to C$ under $G$.
In other words, we have a cartesian square
\[
\xymatrix{
Y \ar[r]^f \ar[d]^b & X \ar[d]^\beta \\
\calC \ar[r]^h & C. \\
}
\]
For any $\sigma \in Z^1(k,G)$ with $\calC^\sigma(k) \ne \emptyset$,
the twisted morphism $b^\sigma \colon Y^\sigma \to \calC^\sigma$
is just like $\beta\colon X \to C$,
since in Section~\ref{S:construction}
we could have replaced $\gamma$ with the composition
$D \injects \calC^\sigma \stackrel{h^\sigma}\to C \stackrel{\gamma}\to \PP^1$
for any connected component $D$ of $\calC^\sigma$ containing a $k$-point;
thus $Y^\sigma(\Adeles)^{\Br}$ 
contains $\calV_\infty(\Adeles) \times \calC^\sigma(k)$
and $f^\sigma(Y^\sigma(\Adeles)^{\Br})$ contains
$\calV_\infty(\Adeles) \times h(\calC^\sigma(k))$.
Taking the union over all such $\sigma$,
and applying the analogue of~\eqref{E:coverings} for $h \colon \calC \to C$,
we see that 
\[
        \Union_{[\sigma] \in H^1(k,G)} f^\sigma(Y^\sigma(\Adeles)^{\Br})
\]
contains $\calV_\infty(\Adeles) \times C(k)$.
Finally, intersect over all $G$ and all $f\colon Y \to X$.
\end{proof}

\begin{remark}
Suppose that $k$ is a number field.
As mentioned in Section~\ref{S:comparisons},
Demarche has announced a proof that 
$X(\Adeles)^{\et,\Br} \subseteq X(\Adeles)^{\descent}$
holds for every nice $k$-variety $X$ \cite{Demarche2008-preprint}.
Assuming this, Theorem~\ref{T:no descent obstruction}
implies that $X(\Adeles)^{\descent}$ is nonempty for our $X$,
and in particular that even the descent obstruction is insufficient
to explain all failures of the Hasse principle.
\end{remark}

\begin{remark}
It is not true that $\beta$ induces an isomorphism $H^1(C,G) \to H^1(X,G)$
for every linear algebraic group $G$.
It fails for $G=\G_m$, for instance, as pointed out to me 
by Colliot-Th\'el\`ene:
the composition $\Pic C \to \Pic X \to \Pic X_\eta \isom \Z$ is zero
but $\Pic X \to \Pic X_\eta$ is nonzero.

So the proof of Theorem~\ref{T:no descent obstruction}
does not directly generalize to prove $X(\Adeles)^{\descent} \ne \emptyset$.
\end{remark}

\begin{remark}
In~\cite{Colliot-Thelene1997} it is conjectured that for every nice variety
over a number field, the Brauer-Manin obstruction is the only 
obstruction to the existence of a zero-cycle of degree~$1$.
In response to an early draft of this paper, Colliot-Th\'el\`ene
has verified this conjecture for the $3$-folds $X$ 
we constructed~\cite{Colliot-Thelene2008-preprint}*{Th\'eor\`eme~3.1}.
\end{remark}

\appendix
\section{Group cohomology}\label{A:cohomology}

In this section, $G$ is a profinite group and
\begin{equation}
\label{E:2-extension}
        0 \to \Z \to A \stackrel{\phi}\to B \to \Z \to 0
\end{equation}
is an exact sequence of discrete $G$-modules, 
with $G$ acting trivially on each copy of $\Z$.
Let $C=\phi(A)$.
Split \eqref{E:2-extension} into short exact sequences
\begin{align*}
        0 \to \Z \to A &\to C \to 0 \\
                     0 &\to C \to B \to \Z \to 0.
\end{align*}
Take cohomology and use $H^1(G,\Z)=\Hom_{\conts}(G,\Z)=0$ 
to obtain exact sequences
\begin{align}
\nonumber  0 \to H^1(G,A) &\to H^1(G,C) \stackrel{\delta_2}\to H^2(G,\Z) \\
\label{E:delta_1}
           H^0(G,B) \to \Z &\stackrel{\delta_1}\to H^1(G,C) \to H^1(G,B) \to 0
\end{align}
Let $\xi = \delta_2(\delta_1(1)) \in H^2(G,\Z)$.
(Thus $\xi$ is the class of the $2$-extension~\eqref{E:2-extension}.)

\begin{lemma}
\label{L:nonzero H^2}
If $H^0(G,B) \to \Z$ is not surjective and $H^1(G,A)=0$, then $\xi \ne 0$.
\end{lemma}

\begin{proof}
Non-surjectivity of $H^0(G,B) \to \Z$ implies $\delta_1(1) \ne 0$,
and $H^1(G,A)=0$ implies that $\delta_2$ is injective.
Thus $\xi \ne 0$.
\end{proof}

\begin{lemma}
\label{L:2-extension}
If $\xi \ne 0$, and the injection $\Z \to A$ and the surjection $B \to \Z$
are both split after restriction to an index-$2$ open subgroup $H$ of $G$,
then the homomorphism $H^1(G,A) \to H^1(G,B)$ induced by $\phi$
is an isomorphism.
\end{lemma}

\begin{proof}
Let $I$ be the image of $\delta_2$
in $H^2(G,\Z) \isom \Hom_{\conts}(G,\Q/\Z)$.
Let $J$ be the image of $\delta_1$.
The splitting of the injection implies that 
$I \subseteq \ker\left(\Hom_{\conts}(G,\Q/\Z) \to \Hom_{\conts}(H,\Q/\Z)\right)
\isom \Hom(G/H,\Q/\Z)$, so $\#I \le 2$.
The splitting of the surjection implies that $2J=0$, 
but $J$ is cyclic, so $\# J \le 2$.
Since $\xi \ne 0$, the composition
$\Z \stackrel{\delta_1}\to H^1(G,C) \stackrel{\delta_2}\to H^2(G,\Z)$
is nonzero, so the induced map $J \to I$ is nonzero.
Therefore $\#I=\#J=2$ and $J \to I$ is an isomorphism.
In particular, $H^1(G,C) \isom H^1(G,A) \directsum J$,
and \eqref{E:delta_1} then yields
\[
        0 \to J \to H^1(G,A) \directsum J \to H^1(G,B) \to 0,
\]
with $J$ mapping identically to $0 \directsum J$.
Thus the map $H^1(G,A) \to H^1(G,B)$ is an isomorphism.
\end{proof}

\section*{Acknowledgements} 

I thank Jean-Louis Colliot-Th\'el\`ene for several comments,
and in particular for suggesting the algebraic
approach in Section~\ref{S:etale covers}
replacing my original topological argument using 
fundamental groups of complex varieties, 
which worked only in characteristic~$0$.
I thank also Martin Olsson.

\end{document}